\documentclass[letterpaper, 10 pt, conference]{IEEEtran}  

\IEEEoverridecommandlockouts                              

\usepackage{epsfig} 
\usepackage{epstopdf}
\usepackage{amsmath} 
\usepackage{amssymb}  
\usepackage{caption}

\usepackage[caption=false,font=footnotesize]{subfig}

\DeclareCaptionLabelSeparator{periodspace}{.\quad}
\newcommand{\squeezeup}{\vspace{-2.5mm}}

\title{\LARGE \bf
Attitude and Angular Velocity Tracking for a Rigid Body using Geometric Methods on the Two-Sphere}

\author{Michalis Ramp$^{1}$ and Evangelos Papadopoulos$^{2}$ 
\thanks{*Financial support by the European Union (European Social Fund-ESF) and Greek national funds through the Operational Program ''Education and Lifelong Learning'' of  the National Strategic Reference Framework Research Funding Program: THALES: Reinforcement of the interdisciplinary and/or interinstitutional research $\&$ innovation is acknowledged.}
\thanks{$^{1}$M. Ramp is with the Department of Mechanical Engineering, National Technical University of Athens, (NTUA) 15780 Athens, Greece
        {\tt\small rampmich@mail.ntua.gr}}%
\thanks{$^{2}$E. Papadopoulos is with the Department of Mechanical Engineering, NTUA, 15780 Athens (tel: +30-210-772-1440; fax: +30-210-772-1455)
        {\tt\small egpapado@central.ntua.gr}}%
}

\begin{document}


\IEEEpubid{\copyright~2015 IEEE. Personal use is permitted. For any other purposes, permission must be obtained from the IEEE. DOI: 10.1109/ECC.2015.7331033}
\IEEEpubidadjcol
\maketitle
\pagestyle{empty}

\begin{abstract}
The control task of tracking a reference pointing direction (the attitude about the pointing direction is irrelevant) while obtaining a desired angular velocity (PDAV) around the pointing direction using geometric techniques is addressed here.
Existing geometric controllers developed on the two-sphere only address the tracking of a reference pointing direction while driving the angular velocity about the pointing direction to zero.
In this paper a tracking controller on the two-sphere, able to address the PDAV control task, is developed globally in a geometric frame work, to avoid problems related to other attitude representations such as unwinding (quaternions) or singularities (Euler angles).
An attitude error function is constructed resulting in a control system with desired tracking performance for rotational maneuvers with large initial attitude/angular velocity errors and the ability to negotiate bounded modeling inaccuracies.
The tracking ability of the developed control system is evaluated by comparing its performance with an existing geometric controller on the two-sphere and by numerical simulations, showing improved performance for large initial attitude errors, smooth transitions between desired angular velocities and the ability to negotiate bounded modeling inaccuracies.
\end{abstract}

\section{INTRODUCTION}
The task of tracking a reference pointing direction (the attitude about the pointing direction is irrelevant) while obtaining a desired angular velocity (PDAV) around the pointing direction represents a fundamental control problem for a variety of robotic applications.
Tiltrotor aircrafts combine the functionality of a conventional helicopter with the long-range, high-velocity performance of a turboprop airplane due to the tilting of the propellers while maintaining a desired propeller speed \cite{Osprey}.
A plethora of military applications involve orienting surveillance apparatus (radar, sonar, sensors) while maintaining a reference angular velocity around the pointing direction in order to scan an area or achieve other objectives gaining a tactical advantage.
In space, a Passive Thermal Control (PTC) technique, also known as "barbecue" roll, is employed in which as a spacecraft is pointed towards a desired direction, it rotates also about an axis to ensure even heat distribution across its surface (the surfaces in direct sunlight reach 390$^o$F while those in the shade -144$^o$F) \cite{PCT}.

Prior work on PDAV is based on conventional representations of attitude like quaternions and Euler angles.
Quaternions exhibit ambiguities in attitude representations since the special orthogonal group SO(3) is double covered, meaning that a global attitude is described by two antipodal points on the three-sphere S$^3$ \cite{EF}, giving rise to unwinding phenomena where the rigid body rotates unnecessarily even though its attitude is extremely close to the desired orientation \cite{RBAC}.
Furthermore a quaternion controller is discontinuous when applied.
Euler angles are defined only locally and exhibit kinematic singularities \cite{O'Reilly}.
Dynamic systems that evolve on nonlinear manifolds cannot be described globally with Euclidean spaces \cite{GCBOOK}.
By utilizing geometric control techniques, control systems are developed by inherently including the properties of the systems nonlinear manifolds in the characterization of the configuration manifold to avoid singularities and ambiguities associated with minimal representations of attitude.
This methodology has been applied to fully/under actuated dynamic systems on Lie groups to achieve almost global asymptotic stability \cite{EF},\cite{RBAC},\cite{GCBOOK},\cite{QUAD},\cite{RAS},\cite{HE},\cite{obstruction}.

In this paper the PDAV control task is addressed in a geometric manner.
Existing geometric controllers developed on the two-sphere, $\text{S}^{2}$, only address the tracking of a reference pointing direction while driving the angular velocity about the pointing direction to zero.
A tracking controller on $\text{S}^2$, able to address the PDAV control task, is developed globally in a geometric framework, avoiding problems related to other attitude representations such as unwinding (quaternions) or singularities (Euler angles).
Inspired by \cite{EF} we develop an attitude error function resulting in a control system with desired tracking performance for rotational maneuvers with large initial attitude/angular velocity errors and the ability to negotiate bounded modeling inaccuracies.
The tracking ability of the developed control system is evaluated by comparing its performance with an existing geometric controller on $\text{S}^{2}$ and by numerical simulations, showcasing improved performance for large initial attitude errors, smooth transitions between desired angular velocities and the ability to negotiate bounded modeling inaccuracies.
In the authors best knowledge a geometric PDAV controller is proposed for the first time.
\IEEEpubidadjcol

\section{KINETICS}
The attitude dynamics of a fully actuated rigid body are studied first.
A body fixed frame $\mathbf{I}_{b}\big\{\mathbf{e}_1,\mathbf{e}_2,\mathbf{e}_3\big\}$, located at the center of mass of the rigid body together with an inertial reference frame $\mathbf{I}_{R}\big\{\mathbf{E}_1,\mathbf{E}_2,\mathbf{E}_3\big\}$, are defined.
The configuration of the rigid body is the orientation of a reference principal axis of the body's body-fixed frame with respect to the inertial frame given by,
\squeezeup
\begin{IEEEeqnarray}{rCl}
\label{eq:1}
\mathbf{q}&=&\mathbf{Q}\mathbf{e}_{3}\IEEEeqnarraynumspace
\end{IEEEeqnarray}
The configuration space is S$^2$=$\{ \mathbf{q}\in\mathbb{R}^{3}\lvert\mathbf{q}^{T}\mathbf{q}=1 \}$.
The plane tangent to the unit sphere at $\mathbf{q}$ is the tangent space $\text{T}_q\text{S}^2=\{ \boldsymbol{\xi}\in\mathbb{R}^{3}\lvert\mathbf{q}^{T}\boldsymbol{\xi}=0 \}$.

For the PDAV application, the attitude configuration is sufficiently described by the unit vector $\mathbf{q}$, but the rigid body's full attitude configuration is defined by the rotation matrix $\mathbf{Q}{\in}$SO(3), mapping a configuration vector from $\mathbf{I}_{b}$ to $\mathbf{I}_{R}$.
It's associated angular velocity vector in $\mathbf{I}_{b}$ is given by,
\begin{IEEEeqnarray}{rCl}
\label{eq:2}
({}^{b}\boldsymbol{\omega})^{\times}&=&\mathbf{Q}^{T}\dot{\mathbf{Q}}
\end{IEEEeqnarray}
where the mappings, $(.)^{\times}$ and its inverse $(.)^{\vee}$ are found in the Appendix, see (\ref{r}).
An additional transitive rotation matrix $\mathbf{Q}_{T}{\in}$SO(3) is utilized, mapping a configuration vector in $\mathbf{I}_{R}$ to a different vector again in $\mathbf{I}_{R}$.
Its associated angular velocity in $\mathbf{I}_{R}$ is given by,
\begin{IEEEeqnarray}{rCl}
\label{eq:3}
(\boldsymbol{\omega})^{\times}&=&\dot{\mathbf{Q}}_{T}(\mathbf{Q}_{T})^{T}
\end{IEEEeqnarray}
The following kinematic equation can be shown using (\ref{eq:3}),
\begin{IEEEeqnarray}{rCl}
\label{eq:dq}
\dot{\mathbf{q}}&=&(\boldsymbol{\omega})^{\times}\mathbf{q}
\end{IEEEeqnarray}
where $\boldsymbol{\omega}{=}\mathbf{Q}{}^{b}\boldsymbol{\omega}$.
The equations of motion are given by,
\begin{IEEEeqnarray}{rCl}
\label{eq:5}
\mathbf{J}{}^{b}\dot{\boldsymbol{\omega}}+({}^{b}\boldsymbol{\omega})^{\times}\mathbf{J}{}^{b}\boldsymbol{\omega}&=&{}^{b}\mathbf{u}-c{}^{b}\boldsymbol{\omega}+\boldsymbol{\tau}
\end{IEEEeqnarray}
together with (\ref{eq:2}).
$\mathbf{J}{\in}\mathbb{R}^{3\times3}$ is the diagonal inertial matrix of the body in $\mathbf{I}_b$ while $^b\mathbf{u}{\in}\mathbb{R}^{3}$ is the applied control moment in $\mathbf{I}_b$.
The second term on the right hand side of (\ref{eq:5}) is a dynamic friction attributed moment with $c{>}0$ while $\boldsymbol{\tau}{\in}\mathbb{R}^{3}$ is a moment depending on the application and it's included for completeness (in a tiltrotor application for example it might be a moment due to gravity while in a space application it might be an interaction moment).
\section{Geometric tracking control on S$^2$}
A control system able to follow a desired smooth pointing direction $\mathbf{q}_{d}(t){\in}\text{S}^2$ and angular velocity around $\mathbf{q}_{d}(t)$, $^b\boldsymbol{\omega}_{d}$ is developed next.
The tracking kinematics equations are,
\begin{IEEEeqnarray}{C}
\label{eq:eofm}
\dot{\mathbf{q}}_{d}=(\boldsymbol{\omega}_{d})^{\times}\mathbf{q}_{d}\;,\;\boldsymbol{\omega}_{d}=\mathbf{Q}_{d}{}^b\boldsymbol{\omega}_{d}
\end{IEEEeqnarray}
Onward the subscript $(.)_{d}$ denotes a desired vector/matrix.
\subsection{Error Function}
The error function is the cornerstone of the design procedure of a control system on a manifold with the performance of the controller directly depending on the function \cite{EF},\cite{GCBOOK}.
To construct a control system on $\text{S}^2$ it is crucial to select a smooth positive definite function $\Psi(\mathbf{q},\mathbf{q}_{d}){\in}\mathbb{R}$ that quantifies the error between the current pointing direction and the desired one.
Using $\Psi(\mathbf{q},\mathbf{q}_{d})$, a configuration pointing error vector and a velocity error vector defined in $\text{T}_q\text{S}^2$ are calculated.
Finally by utilizing a Lyapunov candidate written in terms of the error function and the configuration error vectors, the control procedure is similar to nonlinear control design in Euclidean spaces where the control system is meticulously designed through Lyapunov analysis on $\text{S}^2$ \cite{GCBOOK}.

Almost global reduced attitude stabilizing controllers on $\text{S}^2$ have been studied in \cite{RBAC},\cite{RAS},\cite{CONS},\cite{SA},\cite{HE}.
In \cite{RBAC},\cite{RAS} a control system that stabilizes a rigid body to a fixed reference direction $\mathbf{q}_{d}{\in}\text{S}^2$ is summarized, where the error function used in \cite{RAS} is,
\squeezeup
\begin{IEEEeqnarray}{c}
\label{eq:7}
\Psi_{r}(\mathbf{q},\mathbf{q}_{d})=1-\mathbf{q}^{T}\mathbf{q}_{d}
\end{IEEEeqnarray}
However this error function produces a configuration error vector $\mathbf{e}_{r}{=}{(\mathbf{q}_{d})^{\times}\mathbf{q}}\in\mathbb{R}^{3}$, with non proportional magnitude relative to the current pointing attitude $\mathbf{q}$ and the desired one $\mathbf{q}_{d}$ (Fig. 1(a)).
Observing Fig. 1(a), we see that $\mathbf{e}_{r}{=}0$ not only when $\mathbf{q}{=}\mathbf{q}_{d}$ but also at the antipodal point, where $\mathbf{q}{=}{-}\mathbf{q}_{d}$.
As a consequence, the controller performance degrades because the control effort generated due to $\mathbf{e}_{r}$ diminishes as the initial attitude becomes larger than $\pi{/}2$, reducing its effectiveness  for large attitude errors nullifying thus the main advantage of utilizing geometric control methodologies.
Furthermore the control problem in \cite{RBAC},\cite{RAS}, was the attitude stabilization to a fixed pointing direction while driving the angular velocity to zero.
We address the PDAV control task, where non fixed pointing directions/angular velocities are constantly tracked.

Inspired by \cite{EF}, we modify (7) to avoid the above drawback, and to improve the tracking performance for large initial attitude errors.
For a desired pointing direction/angular velocity  tracking command $(\mathbf{q}_{d},{}^{b}\boldsymbol{\omega}_{d})$ and a current orientation/angular velocity $(\mathbf{q},{}^{b}\boldsymbol{\omega})$, we define the new attitude error function as,
\squeezeup
\begin{IEEEeqnarray}{c}
\label{eq:8}
\Psi(\mathbf{q},\mathbf{q}_{d})={2}-\frac{2}{\sqrt{2}}\sqrt{1+\mathbf{q}^{T}\mathbf{q}_{d}}
\end{IEEEeqnarray}
For any $\mathbf{q}{\in}\text{S}^2$ its dot product with $\mathbf{q}_{d}$ is bounded by $-1{\leqq}\mathbf{q}^{T}\mathbf{q}_{d}{\leqq}1$. Thus $\Psi{>}0$ while $\Psi{=}0$ only if $\mathbf{q}{=}\mathbf{q}_{d}$, i.e $\Psi$ is positive definite about $\mathbf{q}{=}\mathbf{q}_{d}$.
For $\mathbf{q}_{d}$ fixed, the left trivialized derivative of $\Psi$ is
\begin{IEEEeqnarray}{rCl}
\label{eq:9}
\text{T}^{*}\text{L}_{q}&=&\frac{1}{\sqrt{2}\sqrt{1+\mathbf{q}^{T}\mathbf{q}_{d}}} (\mathbf{q}_{d})^{\times}\mathbf{q}=\mathbf{Q}{}^{b}\mathbf{e}_{q}
\end{IEEEeqnarray}
where ${}^{b}\mathbf{e}_{q}$ is the attitude error vector.
We calculate (\ref{eq:9}), utilizing the infinitesimal variation of $\mathbf{q}\in\text{S}^2$ using the exponential map (\ref{expm}) as,
\begin{IEEEeqnarray}{c}
\label{eq:10}
\delta \mathbf{q}=\frac{d}{d\epsilon}\Big\lvert_{\epsilon=0}\text{exp}(\epsilon\boldsymbol{\xi}^{\times})\mathbf{q}=(\boldsymbol{\xi})^{\times}\mathbf{q}
\end{IEEEeqnarray}
where $\boldsymbol{\xi}\in\mathbb{R}^{3}$, to get the derivative of $\Psi$ with respect to $\mathbf{q}$,
\begin{IEEEeqnarray}{l}
\label{eq:11}
\mathbf{D}_{q}\Psi\cdot\delta \mathbf{q}=\frac{1}{\sqrt{2}\sqrt{1+\mathbf{q}^{T}\mathbf{q}_{d}}} (\mathbf{q}_{d})^{\times}\mathbf{q}\cdot\boldsymbol{\xi}=\mathbf{Q}{}^{b}\mathbf{e}_{q}\cdot\boldsymbol{\xi}
\end{IEEEeqnarray}
Notice that in (\ref{eq:9}), (\ref{eq:11}) we introduced the rotation matrix $\mathbf{Q}$ but without affecting/changing the result.
This was done because our equations of motion are written in $\mathbf{I}_{b}$.
Thus the attitude error vector ${}^{b}\mathbf{e}_{q}$, emerges as,
\begin{IEEEeqnarray}{l}
\label{eq:12}
{}^{b}\mathbf{e}_{q}(\mathbf{q},\mathbf{q}_{d},\mathbf{Q})=\frac{1}{\sqrt{2}\sqrt{1+\mathbf{q}^{T}\mathbf{q}_{d}}}\mathbf{Q}^{T}(\mathbf{q}_{d})^{\times}\mathbf{q}
\end{IEEEeqnarray}
and it is well defined in the subset $\text{L}_{{2}}{=}\{\mathbf{q}\in$S$^{2}\lvert\Psi(\mathbf{q},\mathbf{q}_{d})<{2}\}$ since $\mathbf{q}^{T}\mathbf{q}_{d}{>}{-}1$.
Onward this analysis is restricted in $\text{L}_{{2}}$.

The critical points of $\Psi$ are the solutions $\mathbf{q}{\in}\text{S}^2$ to the equation $(\mathbf{q}_{d})^{\times}\mathbf{q}{=}0$ and are given by $\mathbf{q}{=}{\pm}\mathbf{q}_{d}$.
Since $-\mathbf{q}_{d}{\notin}\text{L}_{2}$ only one critical point exist namely $\mathbf{q}_{d}{\in}\text{L}_{2}$.
Utilizing,
\begin{IEEEeqnarray}{C}
\label{eq:12-13}
\mathbf{q}_{d}^{T}\mathbf{q}{=}\lVert\mathbf{q}\rVert\lVert\mathbf{q}_{d}\rVert\cos\theta{=}\cos\theta\IEEEnonumber\\
(\mathbf{q}_{d})^{\times}\mathbf{q}{=}{\lVert\mathbf{q}\rVert}{\lVert\mathbf{q}_{d}\rVert}\sin\theta\mathbf{n}{=}\sin\theta\mathbf{n},\;\mathbf{n}{=}{(\mathbf{q}_{d})^{\times}\mathbf{q}}/{\lVert(\mathbf{q}_{d})^{\times}\mathbf{q}\rVert}\IEEEnonumber
\end{IEEEeqnarray}
with $\theta$ being the angle between $\mathbf{q}_{d}$ and $\mathbf{q}$, we obtain,
\begin{IEEEeqnarray}{rCl}
\label{eq:13}
\Psi={2}(1-\cos\frac{\theta}{2})&,&\lVert\mathbf{e}_{q}\rVert^{2}=\frac{\sin^{2}\theta}{2(1+\cos\theta)}\IEEEnonumber\\
\lVert\mathbf{e}_{q}\rVert^{2}\leqq&\Psi&\leqq{2}\lVert\mathbf{e}_{q}\rVert^{2}\IEEEyesnumber
\end{IEEEeqnarray}
showing that $\Psi$ is locally quadratic. 
Using (4) and (8), the time derivative of $\Psi$ for time varying $\mathbf{q}_{d}(t)$, $\boldsymbol{\omega}_{d}(t)$ is,
\begin{IEEEeqnarray}{rCl}
\label{eq:14}
\dot{\Psi}&=&\frac{1}{\sqrt{2}\sqrt{1+\mathbf{q}^{T}\mathbf{q}_{d}}}(\mathbf{q}_{d})^{\times}\mathbf{q}\cdot(\boldsymbol{\omega}-\boldsymbol{\omega}_{d})\IEEEnonumber\\
&=&\frac{1}{\sqrt{2}\sqrt{1+\mathbf{q}^{T}\mathbf{q}_{d}}}(\mathbf{q}_{d})^{\times}\mathbf{q}\cdot\mathbf{Q}{}^{b}\mathbf{e}_{\omega}
\end{IEEEeqnarray}
where ${}^{b}\mathbf{e}_{\omega}$ is the angular velocity tracking error, given by,
\begin{IEEEeqnarray}{l}
\label{eq:15}
^{b}\mathbf{e}_{\omega}({}^{b}\boldsymbol{\omega},{}^{b}\boldsymbol{\omega}_{d},\mathbf{Q},\mathbf{Q}_{d})={}^{b}\boldsymbol{\omega}-\mathbf{Q}^{T}\mathbf{Q}_{d}{}^{b}\boldsymbol{\omega}_{d}
\end{IEEEeqnarray}
Notice that the angular velocity tracking error lies in  $\text{T}_q\text{S}^2$.
Moreover ${}^{b}\mathbf{e}_{q}$ is well defined in $\text{L}_{2}$ with the attitude error magnitude varying proportionally between pointing orientations $\mathbf{q}$, $\mathbf{q}_{d}$ (see Fig. 1(b)).
Resultantly the control effort generated using ${}^{b}\mathbf{e}_{q}$ will also vary proportionally.
This will improve the tracking performance for angles larger than $\pi{/}2$ in relation to (7) \cite{EF}.
\begin{figure}[!h]
\label{Fig 1. Error functions with vectors}
\centering
\subfloat[Error function $\Psi_{r}$ and error vector $\mathbf{e}_{r}$ (3$^{rd}$ component) with respect to an axis angle rotation.]{\includegraphics[width=0.49\columnwidth]{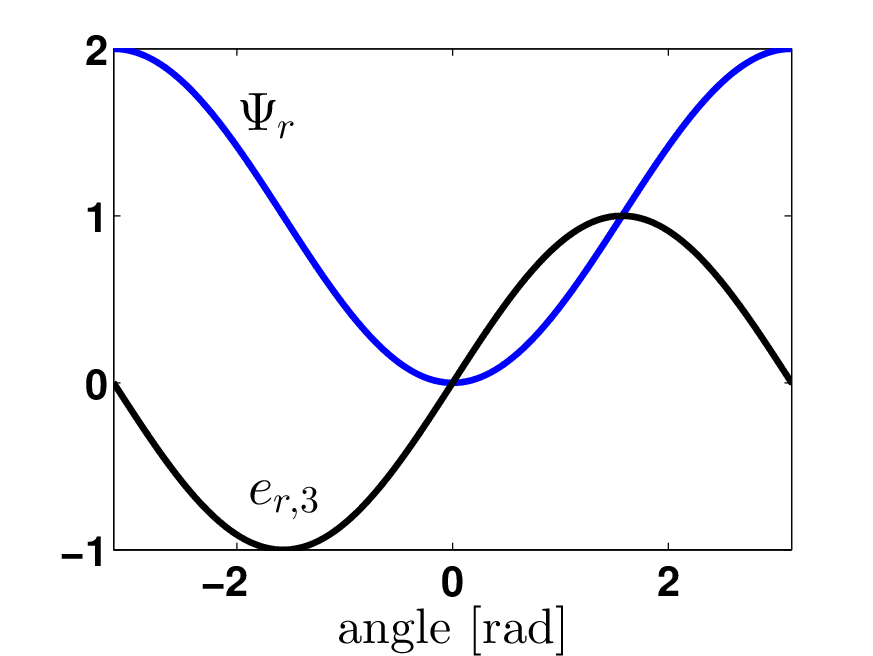}}
~
\subfloat[Derived error function $\Psi$ and error vector ${}^{b}\mathbf{e}_{q}$ (3$^{rd}$ component) with respect to an axis angle rotation.]{\includegraphics[width=0.49\columnwidth]{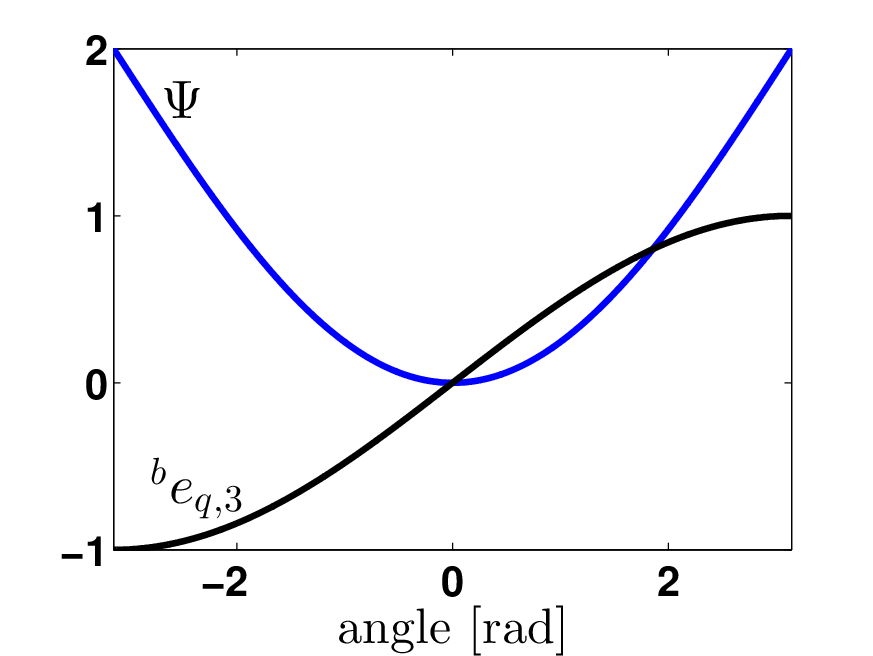}}
\caption{The attitude error function $\Psi_{r}(\mathbf{q},\mathbf{q}_{d})$ together with its error vector $\mathbf{e}_{r}$, (3$^{rd}$ component), from \cite{RAS}, compared with the derived attitude error function $\Psi(\mathbf{q},\mathbf{q}_{d})$, and its error vector ${}^{b}\mathbf{e}_{q}$, (3$^{rd}$ component), as the angle between $\mathbf{q}$ and $\mathbf{q}_{d}$ varies from $-\pi$ to $\pi$.
Only the 3$^{rd}$ component of each error vector is plotted since the others are zero due to the attitude maneuver.}
\end{figure}

\subsection{Error Dynamics}
The error dynamics are developed in order to be used in the subsequent control design.
The derivative of the attitude error function is given in (14).
The derivative of the \textit{configuration error vector} is calculated using (\ref{eq:2}),(\ref{eq:dq}),(\ref{eq:eofm}), after some manipulations,
\begin{IEEEeqnarray}{rCl}
\label{eq:deq}
\frac{d}{dt}{{}^{b}\mathbf{e}_{q}}&{=}&\frac{1}{\sqrt{2}\sqrt{1{+}\mathbf{q}^{T}\mathbf{q}_{d}}}{\mathbf{Q}^{T}}\!\!({({(\boldsymbol{\omega}_{d})^{\times}}\mathbf{q}_{d})^{\times}}\mathbf{q}{+}(\mathbf{q}_{d})^{\times}((\boldsymbol{\omega})^{\times}\mathbf{q}))\IEEEnonumber\\
&&-\frac{1}{2(1{+}\mathbf{q}^{T}\mathbf{q}_{d})}({({(\boldsymbol{\omega}_{d})^{\times}}\mathbf{q}_{d})^{T}}\mathbf{q}{+}\mathbf{q}_{d}^{T}((\boldsymbol{\omega})^{\times}\mathbf{q})){}^{b}\mathbf{e}_{q}\IEEEnonumber\\
&&-({}^{b}\boldsymbol{\omega})^{\times}{}^{b}\mathbf{e}_{q}
\end{IEEEeqnarray}
Using (5), together with $\dot{\mathbf{Q}}_{d}{=}\mathbf{Q}_{d}({}^{b}\boldsymbol{\omega}_{d})^{\times}$ the derivative  of the \textit{angular velocity error} vector is,
\begin{IEEEeqnarray}{rCl}
\label{eq:de_omega}
\frac{d}{dt}{}^{b}\mathbf{e}_{\omega}&=&{}^{b}\dot{\boldsymbol{\omega}}+({}^{b}\boldsymbol{\omega})^{\times}\mathbf{Q}^{T}\mathbf{Q}_{d}{}^{b}\boldsymbol{\omega}_{d}-\mathbf{Q}^{T}\mathbf{Q}_{d}{}^{b}\dot{\boldsymbol{\omega}}_{d}\IEEEnonumber\\
&=&\mathbf{J}^{-1}({}^{b}\mathbf{u}+(\mathbf{J}{}^{b}\boldsymbol{\omega})^{\times}({}^{b}\boldsymbol{\omega})-c{}^{b}\boldsymbol{\omega}+\boldsymbol{\tau})\IEEEnonumber\\
&&+({}^{b}\boldsymbol{\omega})^{\times}\mathbf{Q}^{T}\mathbf{Q}_{d}{}^{b}\boldsymbol{\omega}_{d}-\mathbf{Q}^{T}\mathbf{Q}_{d}{}^{b}\dot{\boldsymbol{\omega}}_{d}
\end{IEEEeqnarray}
\subsection{Pointing direction and angular velocity tracking}
A control system is defined in $\text{L}_{2}$ under the assumption that we have  a fairly accurate estimate of the model parameters showing exponential convergence in an envelope around the zero equilibrium of ${}^{b}\mathbf{e}_{q}$, ${}^{b}\mathbf{e}_{\omega}$ using Lyapunov analysis.

\textit{Proposition 1:} For a desired pointing direction curve $\mathbf{q}_{d}(t)\in\text{S}^{2}$ and a desired angular velocity profile $^{b}\boldsymbol{\omega}_{d}(t)$, around the $\mathbf{q}_{d}(t)$ axis, we define the control moment $^{b}\mathbf{u}$ as,
\begin{IEEEeqnarray}{rCl}
\label{eq:u}
^{b}\mathbf{u}&{=}&\eta^{-1}\widehat{\mathbf{J}}(-\eta(\widehat{\mathbf{f}}+\mathbf{d}){-}(\Lambda+\Psi){}^{b}\dot{\mathbf{e}}_{q}{-}\dot{\Psi}{}^{b}{\mathbf{e}}_{q}{-}\gamma\mathbf{s})\IEEEyessubnumber\\
\mathbf{d}&{=}&({}^{b}\boldsymbol{\omega})^{\times}\mathbf{Q}^{T}\mathbf{Q}_{d}{}^{b}\boldsymbol{\omega}_{d}-\mathbf{Q}^{T}\mathbf{Q}_{d}{}^{b}\dot{\boldsymbol{\omega}}_{d}\IEEEyessubnumber\\
\mathbf{f}&{=}&\mathbf{J}^{-1}((\mathbf{J}{}^{b}\boldsymbol{\omega})^{\times}({}^{b}\boldsymbol{\omega})-c{}^{b}\boldsymbol{\omega}+\boldsymbol{\tau})\IEEEyessubnumber\\
\mathbf{s}&{=}&(\Lambda+\Psi){}^{b}\mathbf{e}_{q}+\eta{}^b\mathbf{e}_{\boldsymbol{\omega}}\IEEEyessubnumber
\end{IEEEeqnarray}
where $\Lambda,\gamma,\eta>0$ positive constants while $\widehat{(.)}$ signifies estimated parameters due to parameter identification errors. 
It will be shown that the above control law stabilizes and maintains ${}^{b}\mathbf{e}_{q}$, ${}^{b}\mathbf{e}_{\omega}$ in a bounded set around the zero equilibrium.
Furthermore for perfect knowledge of the system parameters the above law stabilizes  ${}^{b}\mathbf{e}_{q}$, ${}^{b}\mathbf{e}_{\omega}$ to zero exponentially.

\textit{Proof:} 
See
\cite{Ramp}.

\section{Benchmark Geometric Controller}
A geometric controller from \cite{RAS}, that was derived from (\ref{eq:7}), will be used as a benchmark to gauge the advantages gained from developing the new attitude error function (\ref{eq:8}) and error vector (\ref{eq:12}).

\subsection{Geometric controller on $\text{S}^{2}$}
A stabilization controller, proposed in \cite{RAS}, developed using (\ref{eq:7}) is summarized below,
\begin{IEEEeqnarray}{rCl}
\label{eq:reduced}
^{b}\mathbf{u}&{=}&\mathbf{Q}^{T}({-}K_{r}\mathbf{e}_{r}{-}K_{\omega}\boldsymbol{\omega}){+}({}^{b}\boldsymbol{\omega})^{\times}\mathbf{J}{}^{b}\boldsymbol{\omega}{+}c{}^{b}\boldsymbol{\omega}{-}\boldsymbol{\tau}\IEEEyessubnumber\\
\mathbf{e}_{r}&{=}&(\mathbf{q}_{d})^{\times}\mathbf{q}\IEEEyessubnumber
\end{IEEEeqnarray}
where $K_{r},K_{\omega}{>}0$ are positive constants.
The controller (\ref{eq:reduced}a) has the ability to stabilize asymptotically a rigid body to a fixed pointing direction $\mathbf{q}_{d}{\in}\text{S}^2$, while driving the angular velocity to zero \cite{RAS}.
Furthermore the controller in (\ref{eq:reduced}a), has two additional terms in relation to the actual controller in \cite{RAS} to compensate for the additional moment $\boldsymbol{\tau}$ and the friction type moment $-c{}^{b}\boldsymbol{\omega}$ that exist in the dynamics considered here.
Finally it is written in a more general form than in \cite{RAS}, because there it was specifically developed to stabilize a spherical pendulum.

To have a clear picture of the improvements gained from using the developed error function, the control law (\ref{eq:reduced}a) will not be compared  with (\ref{eq:u}a) but with a controller of similar structure to (\ref{eq:reduced}a).
We replace $\mathbf{e}_{r}$ in (\ref{eq:reduced}a) with $\mathbf{Q}{}^{b}\mathbf{e}_{q}$ to get a similar stabilizing law to (\ref{eq:reduced}a) but with the proposed tracking error (\ref{eq:12}) to get,
\begin{IEEEeqnarray}{l}
\label{eq:29}
{^{b}\mathbf{u}}{=}{\mathbf{Q}^{T}}({-}K_{r}\mathbf{Q}{}^{b}\mathbf{e}_{q}{-}K_{\omega}\boldsymbol{\omega}){+}({}^{b}\boldsymbol{\omega})^{\times}\mathbf{J}{}^{b}\boldsymbol{\omega}{+}c{}^{b}\boldsymbol{\omega}{-}\boldsymbol{\tau}\IEEEyesnumber
\end{IEEEeqnarray}
We do this to compare the error function (\ref{eq:7}), tracking error (\ref{eq:reduced}b) with the ones defined here namely (\ref{eq:8}), (\ref{eq:12}) on equal terms.
The controllers are compared next.
\section{Numerical Simulations}
To highlight the advantages gained from developing the new attitude error function (\ref{eq:8}) and error vector (\ref{eq:12}) we compare (\ref{eq:reduced}a) to (\ref{eq:29}) in a simple stabilization maneuver.

To demonstrate the effectiveness of the proposed controller (\ref{eq:u}a)
we perform a complex PDAV maneuver, which is the main goal of this work, firstly by using the actual model parameters $\mathbf{J}$, $c$, followed by a simulation were (\ref{eq:u}a) uses the estimates 
$\widehat{\mathbf{J}}$, $\widehat{c}$.
This was done to showcase the robustness and the enhanced performance of the proposed controller.
The system parameters are:
\begin{IEEEeqnarray*}{C}
\label{PARAMETERS}
c=0.3\:\text{Nm(s/rad)},\;\widehat{c}=c+0.03c
\IEEEnonumber\\
\mathbf{J}=\text{diag}{\left(
0.0294,0.0305,0.0495
\right)}\:\text{kgm}^2,\;\widehat{\mathbf{J}}=\mathbf{J}+0.14\mathbf{J}
\IEEEnonumber\\
\boldsymbol{\tau}=\mathbf{0}\:\text{Nm}
\end{IEEEeqnarray*}
The controller parameters are chosen through pole placement by choosing time constants/damping coefficients,
\begin{IEEEeqnarray}{C}
\label{30abcde}
K_{r}=\text{diag}{\left(
4.234,4.392,7.128
\right)}\IEEEyessubnumber\\
K_{\omega}=\text{diag}{\left(
7.056,7.320,11.88
\right)}\cdot 10^{-1}\IEEEyessubnumber\\
\Lambda=144,\:\eta=24,\:\gamma=10\IEEEyessubnumber
\end{IEEEeqnarray}
\subsection{Pointing direction stabilization}
We move forward with the comparison of (\ref{eq:reduced}a) and (\ref{eq:29}) to get a clear picture of the improvements gained from the developed error function/vector.
This comparison will take place under the assumption of perfect knowledge of the system parameters.
The initial configuration/conditions are,
\begin{IEEEeqnarray}{C}
\label{IC's}
\mathbf{q}(t{=}0){=}[
0{;}0{;}1],\:\mathbf{Q}(t{=}0){=}\mathbf{I},\:\boldsymbol{\omega}(t{=}0){=}[
0{;}0.3{;}0]\IEEEnonumber\IEEEeqnarraynumspace
\end{IEEEeqnarray}
The controllers must initially stabilize the $\mathbf{e}_{3}$ body fixed axis of the body to the following equilibrium,
\begin{IEEEeqnarray}{C}
\label{REDUCED ATTITUDE STABILIZATION q1}
\mathbf{q}_{d}{=}[
0{;}-0.0175{;}-0.9998],\:\boldsymbol{\omega}{=}[
0{;}0{;}0]\IEEEnonumber\IEEEeqnarraynumspace
\end{IEEEeqnarray}
where the vector $\mathbf{q}_{d}$ denotes a $\text{179}^{o}$ rotation around the $\mathbf{E}_{1}$ axis.
Then the $\mathbf{e}_{3}$ body fixed axis must rotate back $\text{89}^{o}$ around the $\mathbf{E}_{1}$ axis to,
\begin{IEEEeqnarray}{C}
\label{REDUCED ATTITUDE STABILIZATION q2}
\mathbf{q}_{d}{=}[
0{;}-1{;}0],\:\boldsymbol{\omega}{=}[
0{;}0{;}0]\IEEEnonumber\IEEEeqnarraynumspace
\end{IEEEeqnarray}
Examining Fig. 2(a), the effectiveness of (\ref{eq:29}) for large initial attitude errors is demonstrated, as $\Psi$ converges faster to zero.
Also it is observed that the developed error vector (\ref{eq:12}) (Fig. 2(b)) (solid black line) steers the rigid body, Fig. 2(a), and angular velocity, Fig. 2(c), to the desired equilibrium faster.
This can be further substantiated by examining the magnitude of the error vector (black, solid line), Fig. 2(b), which is larger for large initial angle differences in comparison to the one generated from (\ref{eq:reduced}b) (blue, dashed line), resulting to an actively engaged controller in large angle maneuvers.
For initial angles less than $\pi/2$ however, (\ref{eq:reduced}a) drives the system to the desired equilibrium slightly faster (Fig. 2(a,d)).
Nevertheless the inability of (\ref{eq:reduced}a) to swiftly steer the system to the desired attitude when the initial angle difference is larger than $\pi/2$ makes the developed error vector (\ref{eq:12}) more effective as it guarantees a uniform/homogeneous response.
\begin{figure}[!thpb]
\label{Fig 2. Attitude stabilization and tracking}
\centering
\subfloat[ ]{\includegraphics[width=0.49\columnwidth]{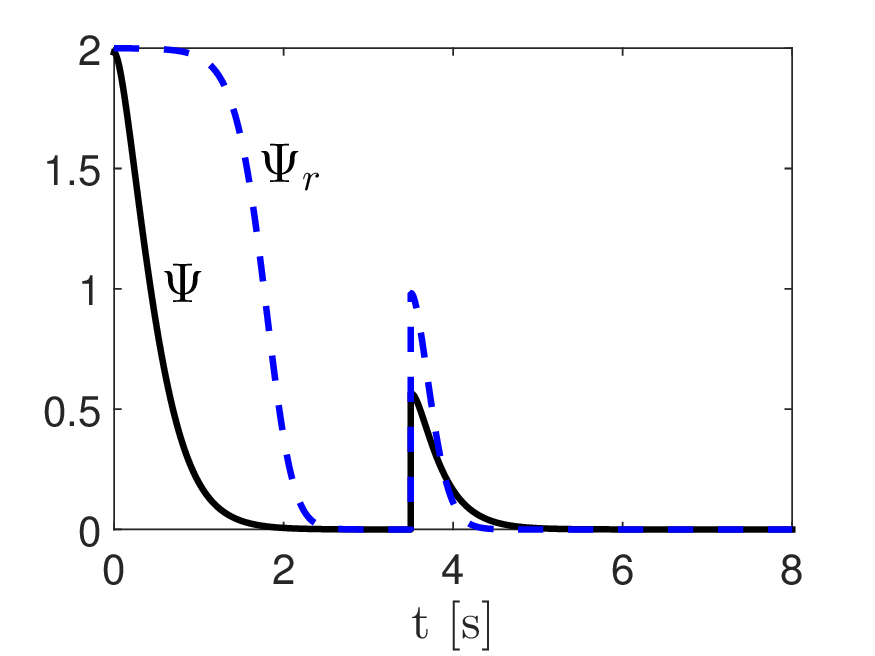}}
~
\subfloat[ ]{\includegraphics[width=0.49\columnwidth]{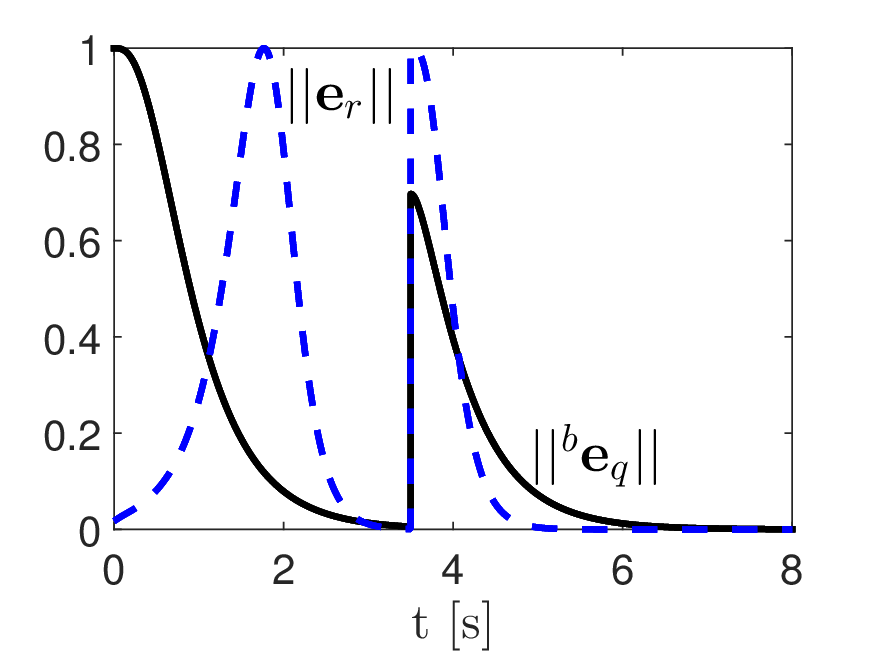}} 

\subfloat[ ]{\includegraphics[width=0.49\columnwidth]{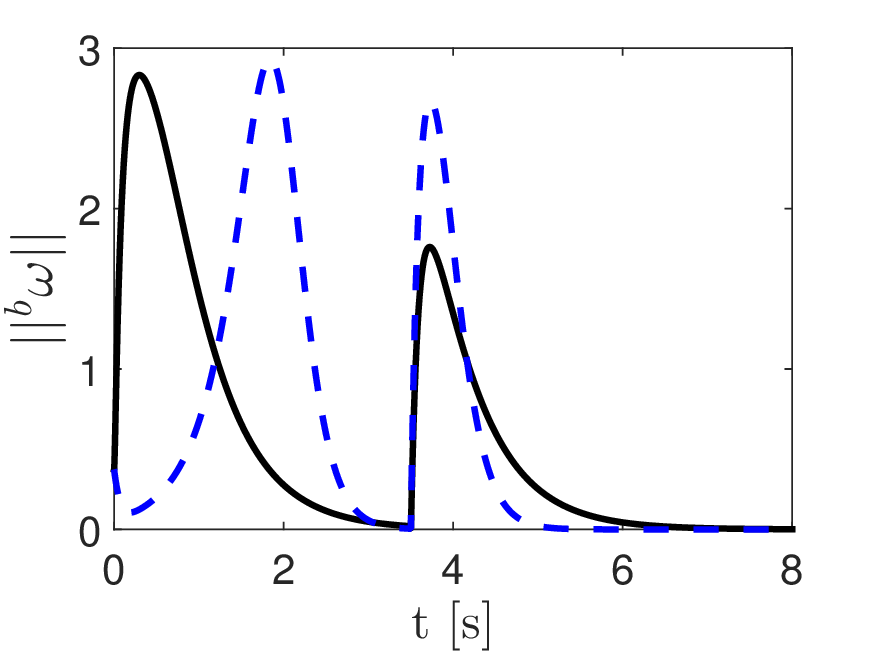}}
~
\subfloat[ ]{\includegraphics[width=0.49\columnwidth]{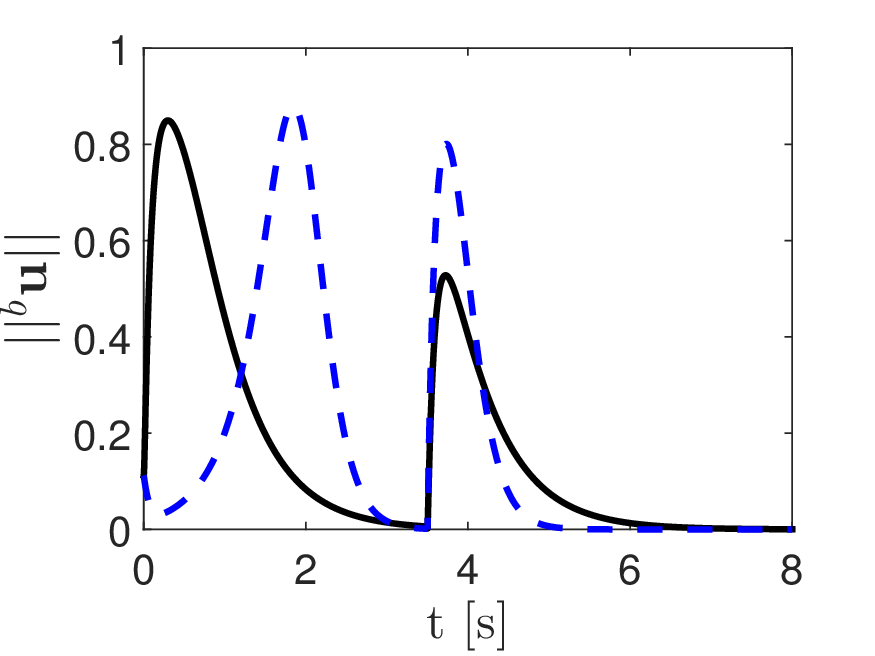}} 
\caption{Pointing direction stabilization where the black solid lines denote the response due to the derived error vector. The blue dashed lines denote the response due to the error vector from \cite{RAS}. Throughout the simulation the angular velocity is driven to zero (a) Attitude comparison through the respective error functions $\Psi_{r}$, $\Psi$. (b) Error vector comparison through the 2-norm. (c) Angular velocity comparison through the 2-norm, (rad/s). (d) Control moment comparison through the 2-norm, (Nm).}
\end{figure}
\subsection{Pointing direction and angular velocity tracking}
A complex PDAV maneuver to test the effectiveness of the proposed controller (\ref{eq:u}a) is performed, firstly by using the actual model parameters $\mathbf{J}$, $c$, followed by a simulation were (\ref{eq:u}a) uses the estimates $\widehat{\mathbf{J}}$, $\widehat{c}$.
The desired attitude trajectory and angular velocity profile about the pointing axis are given in Fig. 3 and generated using smooth polynomials \cite{POLY}, see (\ref{polyn}), with the third Euler angle to be zero and will be omitted from now on (see (\ref{euang})).
The trajectories are,
\begin{IEEEeqnarray*}{C}
\label{PDAV1}
{}^{b}\boldsymbol{\omega}_{d}{=}[
0{;}0{;}{}^{b}\omega_{d,3}(t)],\:\mathbf{Q}_{d}{=}\mathbf{Q}_{313}(\theta(t),\phi(t)),\:\mathbf{q}_{d}{=}\mathbf{Q}_{d}[
0{;}0{;}{}1]
\end{IEEEeqnarray*}

A step command is issued initially where the body fixed axis $\mathbf{e}_{3}$ is required to rotate around the $\mathbf{E}_{1}$ axis 179$^{o}$ (Fig. 3b).
At $t{=}1s$, a smooth trajectory of a $7s$ duration is initialized (see Fig. 3b) where the $\mathbf{e}_{3}$ body fixed axis must rotate back $\text{89}^{o}$ around the $\mathbf{E}_{1}$ axis, while simultaneously the $\mathbf{e}_{3}$ body fixed axis must rotate $\text{90}^{o}$ around the $\mathbf{E}_{3}$ axis (Fig. 3a).
The angular velocity for the first five seconds is smoothly increased up to ${}^{b}\omega_{d,3}{=}10\,\text{rad/s}$ (Fig. 3c).
It is then maintained at that level for five seconds and finally it is driven back to zero.
The desired trajectories contain overlapping parts with both the pointing direction and the angular velocity around it actively regulated.
The PDAV trajectory was chosen in this manner to investigate the ability of the controller in \textit{orienting} the body, while simultaneously regulating the \textit{angular velocity} about $\mathbf{q}_{d}$.

\begin{figure}[!h]
\label{Fig 3. Trajectory and angular velocity profiles}
\centering

\subfloat[ ]{\includegraphics[width=0.32\columnwidth]{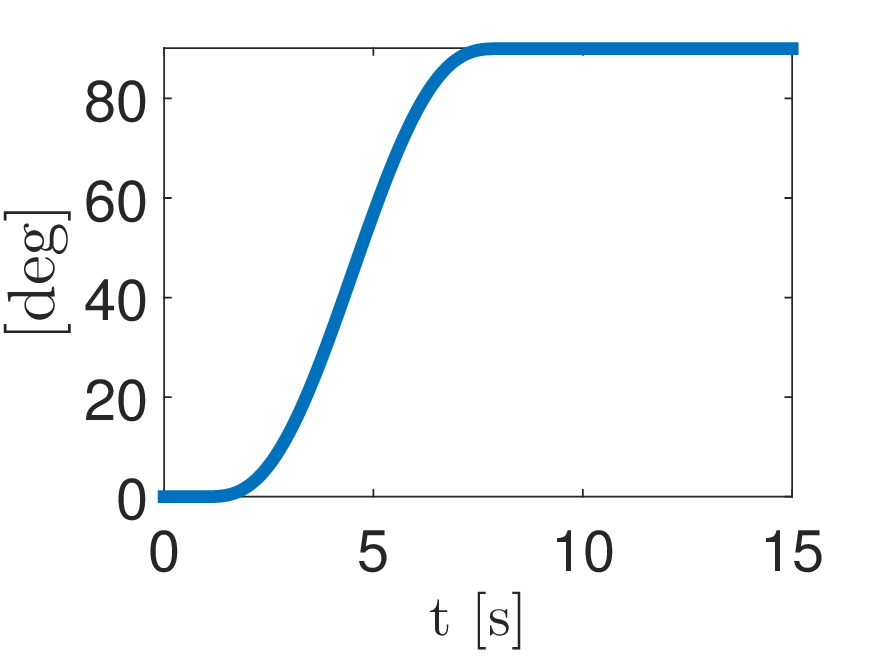}}
~
\subfloat[ ]{\includegraphics[width=0.32\columnwidth]{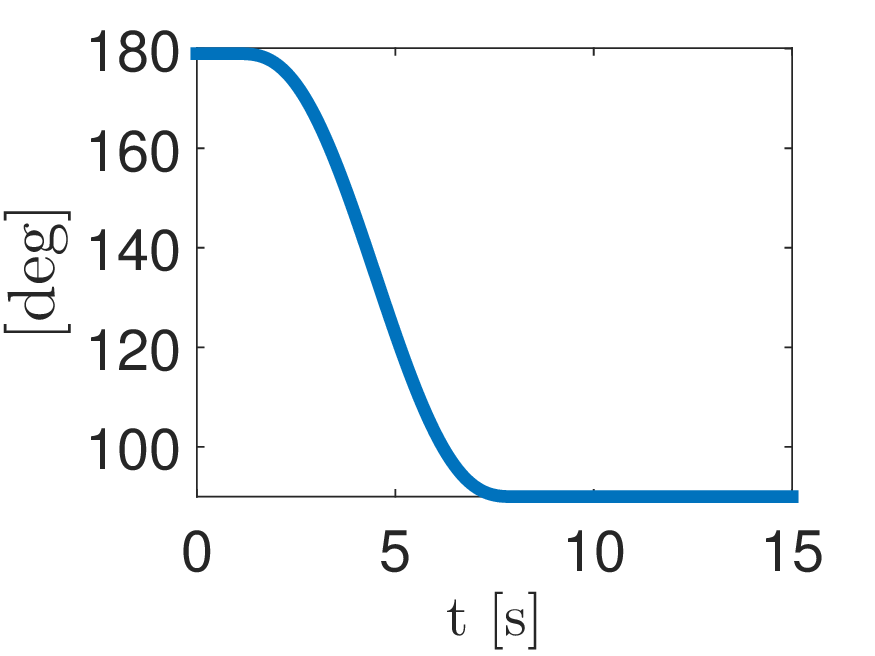}} 
~
\subfloat[ ]{\includegraphics[width=0.32\columnwidth]{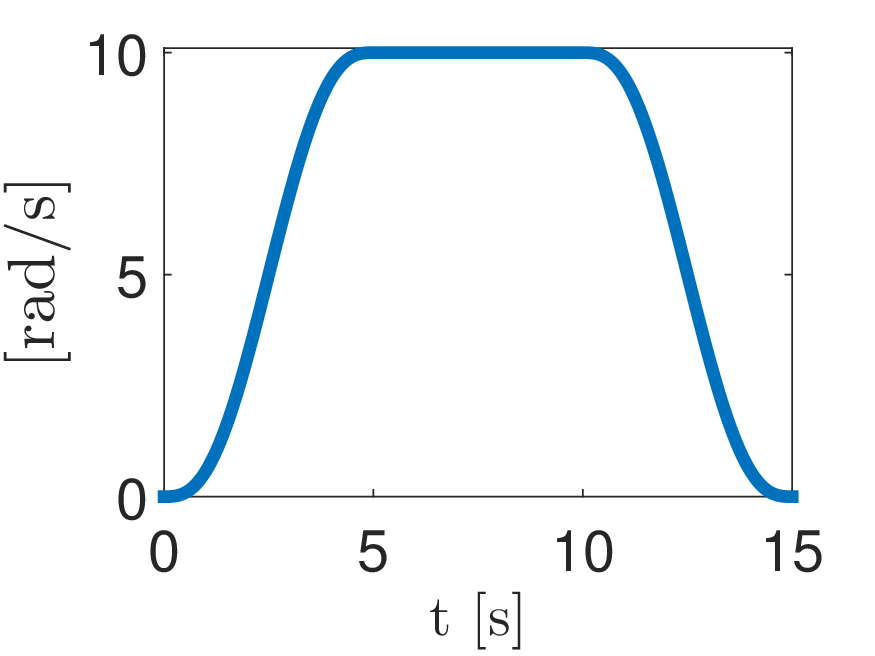}} 
\caption{Rotation matrix design using the 313 Euler angles with $\mathbf{Q}_{d}{=}\mathbf{Q}_{313}(\theta(t),\phi(t))$. The third angle is zero and thus omitted (a) $\phi(t)$. (b) $\theta(t)$. (c) ${}^{b}\omega_{d,3}(t)$. The desired angular velocity is ${}^{b}\boldsymbol{\omega}_{d}{=}[
0{;}0{;}{}^{b}\omega_{d,3}(t)]$.}
\end{figure}

During the first simulation, perfect knowledge of the system is assumed, i.e, $\mathbf{J}$, $c$ are used in (\ref{eq:u}a).
The proposed controller is able to track the desired attitude/angular velocity trajectories very effectively.
Both the pointing direction $\mathbf{q}_{d}$ and the angular velocity about the pointing direction, ${}^{b}\omega_{d,3}$, are tracked almost exactly (Fig. 4(a,b)).
Minor oscillations are observed  during a portion of the tracking maneuver (Fig. 4(b)), specifically when both $\mathbf{q}_{d}$ and ${}^{b}\omega_{d,3}$ are simultaneously varied.
This is due to the desired velocity profile defined by $\mathbf{Q}^{T}\mathbf{Q}_{d}{}^{b}\boldsymbol{\omega}_{d}$ as can be seen in Fig. 4(b) (dashed line) and can be attributed to gyroscopic phenomena.

\begin{figure}[!h]
\label{Fig 4. Pointing direction angular velocity tracking}
\centering
\subfloat[ ]{\includegraphics[width=0.49\columnwidth]{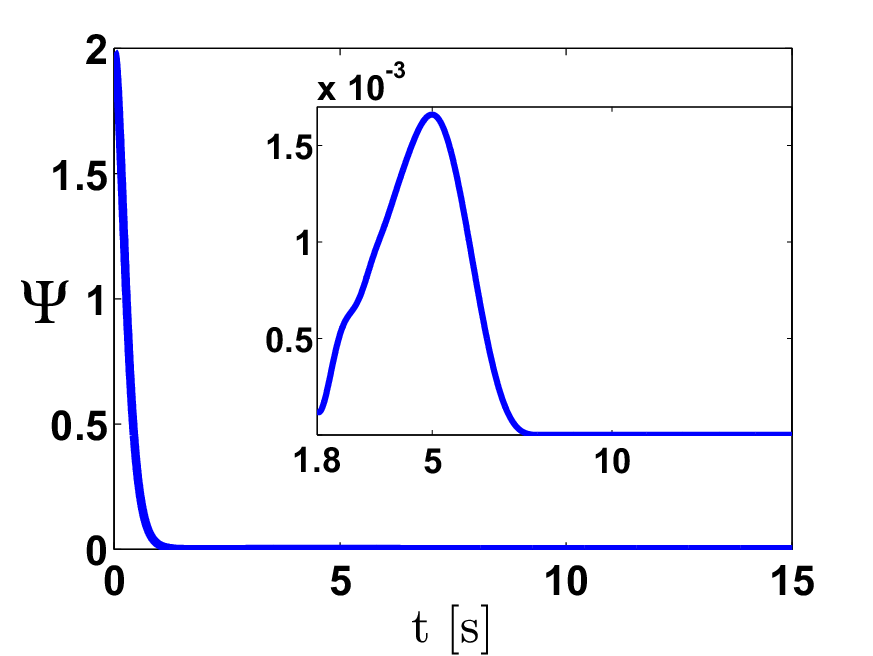}}
~
\subfloat[ ]{\includegraphics[width=0.49\columnwidth]{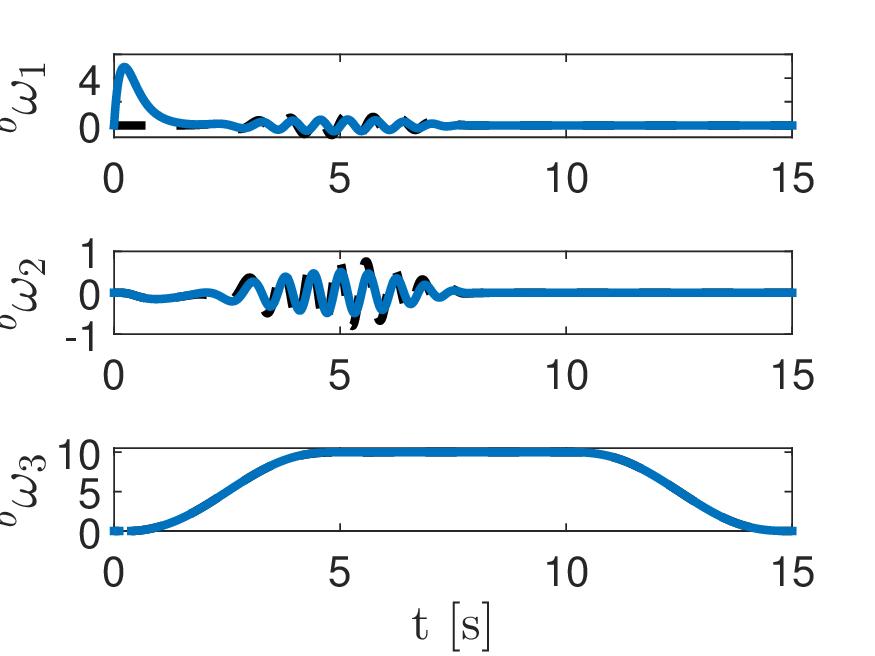}} 

\subfloat[ ]{\includegraphics[width=0.49\columnwidth]{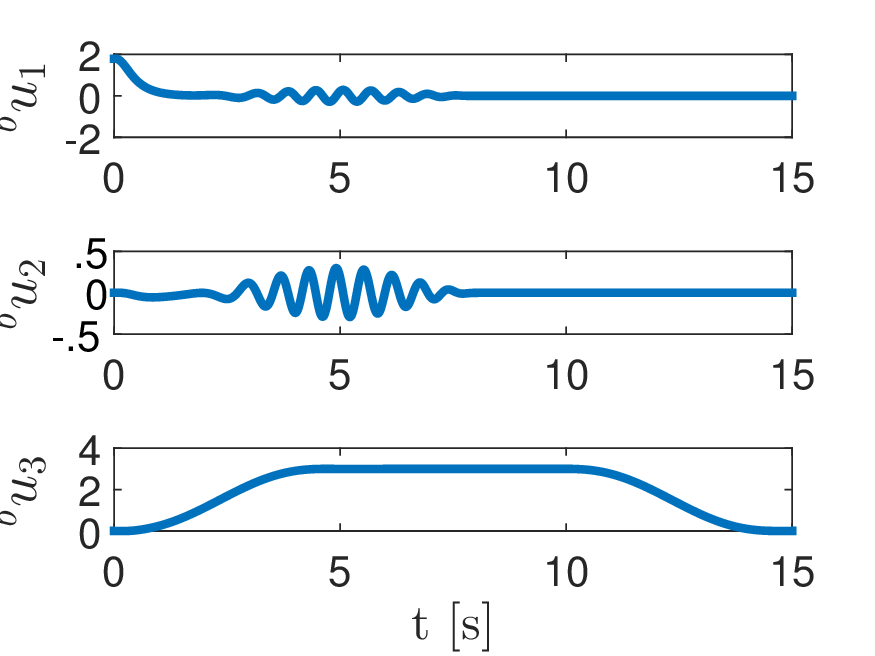}}
~
\subfloat[ ]{\includegraphics[width=0.49\columnwidth]{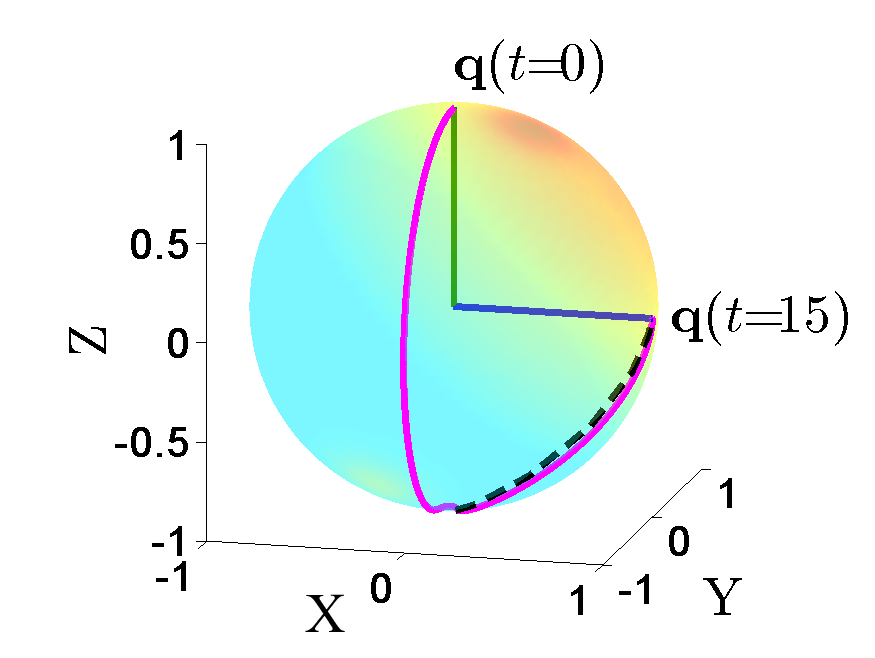}} 
\caption{PDAV tracking using the developed controller (\ref{eq:u}a). Perfect knowledge of the system is assumed. The dashed line indicates desired response (a) Attitude error function response. (b) Angular velocity response, (rad/s). (c) Generated control moment, (Nm). (d) Attitude maneuver on S$^2$.}
\end{figure}

The actual pointing maneuver is indeed smooth and this is apparent in Fig. 4(a) where during the maneuver the attitude error varies smoothly and is maintained below $\Psi{<}1.7{\cdot}10^{-3}$, translating to less than 0.153 deg with respect to an axis angle rotation.
Even though the angular velocity tracking command was designed to be smooth (Fig. 3(c)) the proposed controller doesn't require complicated trajectory commands; it can be shown to transition smoothly from one angular velocity to another using simple step velocity commands.
The claim that the developed controller is smooth, without high frequency chattering, is validated by observing the generated control torque, in Fig. 4(c) which shows smooth control inputs.
The traversed trajectory is shown in Fig. 4(d) demonstrating near perfect tracking of $\mathbf{q}_{d}(t)$, ${}^{b}\boldsymbol{\omega}_{d}(t)$.

To test the robustness of the proposed controller a second simulation experiment is conducted were we repeat the same trajectories as before but the system parameters provided in (\ref{eq:u}a) are estimates $\widehat{\mathbf{J}}$, $\widehat{c}$, with the gains kept the same.
The gains were kept identical to highlight the enhanced performance of the proposed controller but this imposes a bound on how much the estimates used in (\ref{eq:u}a) can be varied.
Despite inaccurate estimates that reach 14$\%$ for some parameters, the proposed controller (\ref{eq:u}a) is able to track the desired attitude/angular velocity trajectories effectively within bounds (Fig. 5(a,b)).
The angular velocity tracking error ${}^{b}\mathbf{e}_{\omega}$ is compared with the one from the previous simulation (Fig. 5(c,d)) revealing that the angular velocity tracking error has increased considerably especially the third component $^{b}e_{\omega,3}$.
This can be remedied by increasing the gain $\gamma$ (depending on available control input bounds).
Finally it can be concluded that the proposed controller can effectively negotiate the PDAV control task in a singularity free-manner, while negotiating bounded parametric inaccuracies exponentially stabilizing ${}^{b}\mathbf{e}_{q}$, ${}^{b}\mathbf{e}_{\omega}$ in a bounded set around the zero equilibrium.
\begin{figure}[!h]
\label{Fig 5. PDAV tracking with J_=0.14, c_=0.03c}
\centering
\subfloat[ ]{\includegraphics[width=0.49\columnwidth]{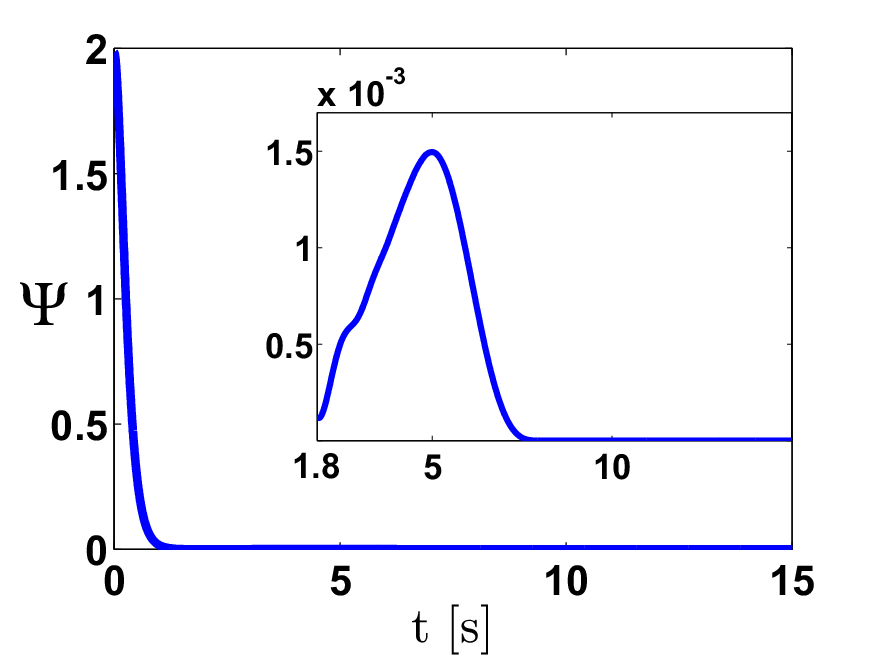}}
~
\subfloat[ ]{\includegraphics[width=0.49\columnwidth]{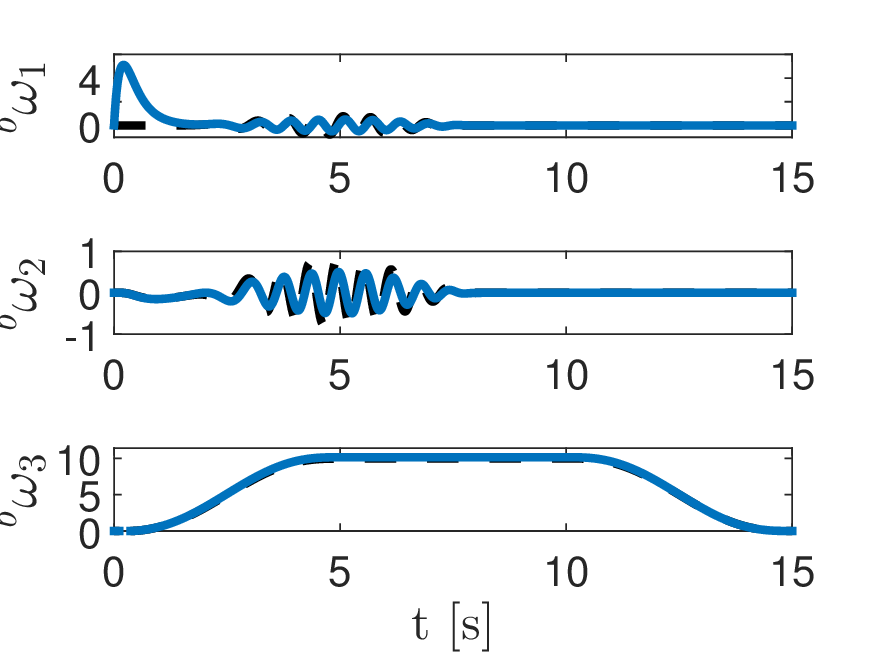}} 

\subfloat[ ]{\includegraphics[width=0.49\columnwidth]{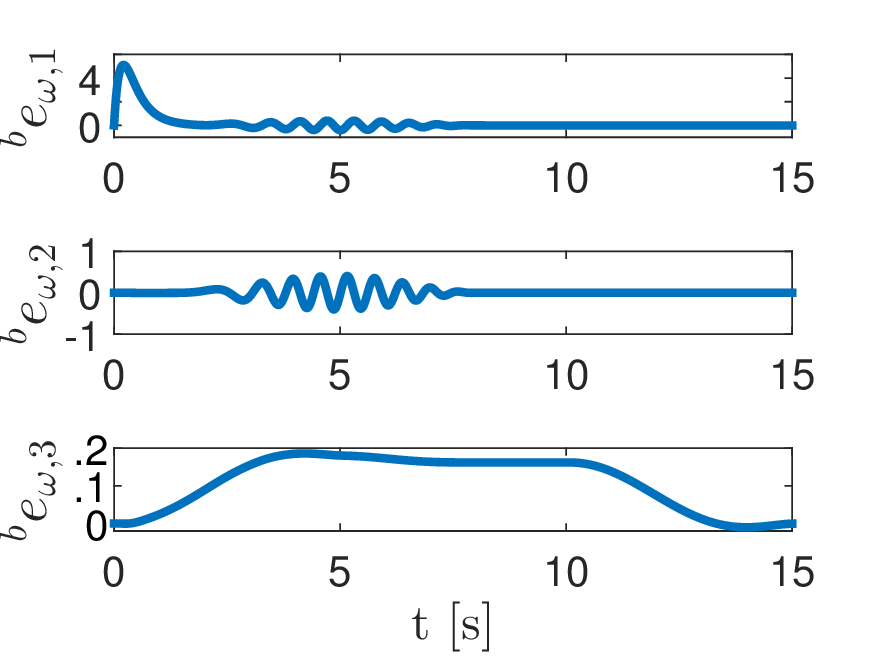}}
~
\subfloat[ ]{\includegraphics[width=0.49\columnwidth]{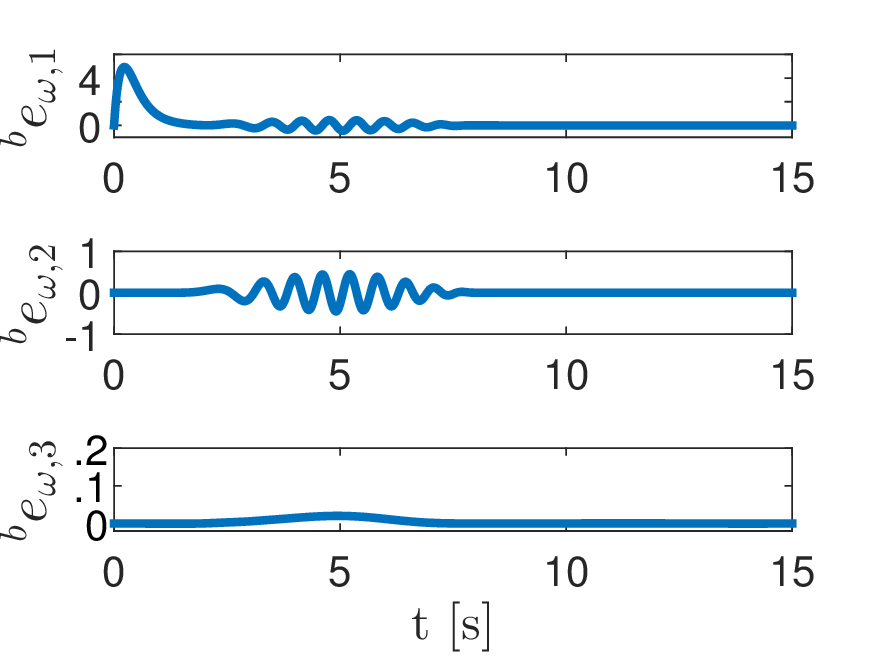}} 
\caption{Simulation using the derived controller (\ref{eq:u}a) but with the system parameters estimates $\widehat{\mathbf{J}}$, $\widehat{c}$, showcasing the ability of the controller to negotiate bounded modeling inaccuracies up to 14$\%$. The dashed line indicates desired response (a) Error function response $\Psi$. The controller utilized $\widehat{\mathbf{J}}$, $\widehat{c}$. (b) Angular velocity response, (rad/s). The controller utilized $\widehat{\mathbf{J}}$, $\widehat{c}$. (c) Angular velocity error vector under the action of (\ref{eq:u}a), (rad/s). The controller utilized $\widehat{\mathbf{J}}$, $\widehat{c}$. (d) Angular velocity error vector under the action of (\ref{eq:u}a), (rad/s). The controller utilized $\mathbf{J}$, $c$.}
\end{figure}
\section{CONCLUSIONS}
In this paper the PDAV control task was addressed in a geometric framework.
An attitude error function was defined and used to develop a singularity-free control system on $\text{S}^{2}$ with improved tracking performance for rotational maneuvers with large initial attitude errors, able to negotiate bounded modeling inaccuracies and exponentially stabilizing ${}^{b}\mathbf{e}_{q}$, ${}^{b}\mathbf{e}_{\omega}$ in a bounded set around the zero equilibrium.
The tracking ability of the developed control system was evaluated by comparing its performance with an existing geometric controller on $\text{S}^{2}$ and by numerical simulations, demonstrating improved tracking performance for large initial orientation errors, smooth transitions between desired angular velocities and the ability to negotiate bounded modeling inaccuracies.




\renewcommand{\theequation}{A\arabic{equation}}
\setcounter{equation}{0}  
\section*{APPENDIX}  
Vector space isomorphism where $\mathbf{r}\in\mathbb{R}^3$
\begin{IEEEeqnarray}{C}
(\mathbf{r})^{\times}{=}[0,-r_3,r_2;r_3,0,-r_1;-r_2,r_1,0],((\mathbf{r})^{\times})^{\vee}{=}\mathbf{r}\IEEEeqnarraynumspace\label{r}
\end{IEEEeqnarray}
Exponential map using the Rodrigues formulation \cite{O'Reilly},
\begin{IEEEeqnarray}{C}
\text{exp}(\epsilon\boldsymbol{\xi}^{\times})=\mathbf{I}+\boldsymbol{\xi}^{\times}\sin{\epsilon}+(\boldsymbol{\xi}^{\times})^{2}(1-\cos{\epsilon})\label{expm}
\end{IEEEeqnarray}
Constants, $[\alpha_{0},{-},\alpha_{5}]$, used to define sixth degree polynomials for trajectory generation \cite{POLY},
\begin{IEEEeqnarray}{L}
{}^{b}\omega_{d,3}(t{\leq}5){=}[0{,}0{,}0{,}0.8{,}{-}0.24{,}0.0192]\label{polyn}\\
{}^{b}\omega_{d,3}(t{\geq}10){=}[5130{,}{-}2160,360{,}{-}29.6{,}1.2{,}{-}0.0192]\IEEEnonumber\\
\phi(1{\leq}t{\leq}8){=}[{-}3.2183{,}10.28{,}{-}11.57{,}5.19{,}{-}0.7229{,}0.0321]\IEEEnonumber\\
\theta(1{\leq}t{\leq}8){=}[182.2{,}{-}10.17{,}11.44{,}{-}5.137{,}0.7149{,}{-}0.0318]\IEEEnonumber
\end{IEEEeqnarray}
Attitude through Euler-Angles ($c\gamma_i=\cos\gamma_i$, $s\gamma_i=\sin\gamma_i$)
\begin{IEEEeqnarray}{C}
\mathbf{Q}_{313}(\gamma_3{=}0{,}\gamma_2{,}\gamma_1){=}
\arraycolsep=1.4pt
{\mathbf{I}}{\cdot}
{\begin{bmatrix}
1&0&0\\
0&c\gamma_2&-s\gamma_2\\
0&s\gamma_2&c\gamma_2
\end{bmatrix}}{\cdot}{\begin{bmatrix}
c\gamma_1&-s\gamma_1&0\\
s\gamma_1&c\gamma_1&0\\
0&0&1
\end{bmatrix}}\label{euang}
\end{IEEEeqnarray}

\end{document}